\documentclass[a4paper,11pt] {amsart}
\usepackage{amssymb,latexsym,amsmath}
\usepackage{amsfonts,amsbsy,bm}
\usepackage{color}
\renewcommand {\phi} {{\varphi}}

\newcommand {\al} {{\alpha}}
\newcommand {\ba} {{\beta}}
\newcommand {\da} {{\delta}}
\newcommand {\Da} {{\Delta}}

\newcommand {\la} {{\lambda}}

\renewcommand {\O} {{\Omega}}

\newcommand {\N} {{\mathbb N}}

\newcommand {\Scirc} {\raise.2ex\hbox{$\scriptstyle\circ$}}

\newcommand {\im} {{\Im\!\mbox{\small\it m}\,}}

\newcommand {\Proof} {\noindent{\bf P{\footnotesize\bf ROOF}: } \ }

\newcommand {\ProofEnd} {
             \begin{flushright} \vskip -0.2in $\Box$ \end{flushright}}

\newcommand{\Ba}[1]{\begin{array}{#1}}
\newcommand{\Ea}{\end{array}}
\newcommand{\Be}{\begin{equation}}
\newcommand{\Ee}{\end{equation}}
\newcommand{\Bea}{\begin{eqnarray}}
\newcommand{\Eea}{\end{eqnarray}}
\newcommand{\Beas}{\begin{eqnarray*}}
\newcommand{\Eeas}{\end{eqnarray*}}
\newcommand{\Benu}{\begin{enumerate}}
\newcommand{\Eenu}{\end{enumerate}}
\newcommand{\Bi}{\begin{itemize}}
\newcommand{\Ei}{\end{itemize}}

\newcommand{\BR}{\begin{Remark} \em}
\newcommand{\ER}{\end{Remark}}
\newcommand{\BE}{\begin{example} \em}
\newcommand{\EE}{\end{example}}
\newcounter{remark}

\newtheorem{theorem}[equation]{T{\hskip 0pt\footnotesize\bf HEOREM}}
\newtheorem{proposition}[equation]{P{\hskip 0pt\footnotesize\bf ROPOSITION}}
\newtheorem{corollary}[equation]{C{\hskip 0pt\footnotesize\bf OROLLARY}}
\newtheorem{lemma}[equation]{L{\hskip 0pt\footnotesize\bf EMMA}}
\newtheorem{Remark}[equation]{R{\hskip 0pt\footnotesize\bf EMARK}}
\newtheorem{definition}[equation]{D{\hskip 0pt\footnotesize\bf EFINITION}}
\newtheorem{example}[equation]{E{\hskip 0pt\footnotesize\bf XAMPLE}}

\newcommand{\bprop} {\begin{proposition}}
\newcommand{\eprop} {\end{proposition}}
\newcommand{\btheo} {\begin{theorem}}
\newcommand{\etheo} {\end{theorem}}
\newcommand{\blem} {\begin{lemma}}
\newcommand{\elem} {\end{lemma}}
\newcommand{\bcor} {\begin{corollary}}
\newcommand{\ecor} {\end{corollary}}

scaled \magstep1

\renewcommand{\O}{\Omega}

\renewcommand{\ba}{{\bf \al}}
\newcommand{{\tBox}}{{\widetilde{\raisebox{-0.2ex}[1.25ex][0ex]{$\Box$}}}}

\begin{document}
\title[Carleson embeddings and two operators on Bergman spaces]{Carleson embeddings and two operators on Bergman spaces of tube
  domains over symmetric cones.\footnotetext{\emph{2000 Math Subject Classification:} 42B35, 32M15.}
\footnotetext{\emph{Keywords}: Bergman projection, Toeplitz operator, Ces\`aro-type operator, Schatten class, Besov space, symmetric
cone.}}
\author{Cyrille Nana}
\address{Cyrille Nana, Department of Mathematics, Faculty of Science, University of Buea, P. O. Box 63, Buea Cameroon}
\email{nana.cyrille@ubuea.cm}
\author{Beno\^it Florent Sehba}
\address{Beno\^it Florent Sehba,  Department of Mathematics, University of Ghana, P. O. Box LG 62 Legon, Accra, Ghana}
\email{bfsehba@ug.edu.gh}
 \maketitle
\begin{abstract}
We prove Carleson embeddings for Bergman spaces of tube domains over symmetric cones, we apply them to characterize symbols of bounded Ces\`aro-type operators from weighted Bergman spaces to weighted Besov spaces. We also obtain Schatten class criteria of Toeplitz operators and Ces\`aro-type operators on  weighted Hilbert-Bergman spaces.
\end{abstract}

\section{Introduction}
\setcounter{equation}{0} \setcounter{footnote}{0}
\setcounter{figure}{0} Given an irreducible symmetric  cone $\Omega$
in $\mathbb R^n$, we denote by $\mathcal D=\mathbb R^n+i\Omega$
the tube domain over $\Omega$. Following the notations of \cite{FK} we
denote by $r$ the rank of the cone $\Omega$ and by $\Delta$
the determinant function of $\mathbb R^n$. The Lorentz cone
$\Lambda_n$ is defined for $n\ge 3$ by
$$\Lambda_n=\{(y_1,\cdots,y_n)\in \mathbb R^n: y_1^2-\cdots-y_n^2>0,\,\,\,y_1>0\};$$
it is a typical example of symmetric cone on $\mathbb R^n$ with rank $2$ and its associated determinant function is given by the
Lorentz form
$$\Delta(y)=y_1^2-\cdots-y_n^2.$$  We shall denote by $\mathcal{H}(\mathcal D)$ the vector space
of holomorphic functions on $\mathcal D$.

\vskip .2cm
 For $1\le p<\infty$ and $\nu \in \mathbb R$, let
$L^p_\nu(\mathcal D)=L^p(\mathcal
D,\Delta^{\nu-\frac{n}{r}}(y)dx\,dy)$
denotes the space of functions $f$ 
satisfying the condition
$$\|f\|_{p,\nu}=||f||_{L^p_\nu(\mathcal D)}:=\left(\int_{\mathcal D}|f(x+iy)|^p\Delta^{\nu-\frac{n}{r}}(y)dxdy\right)^{1/p}<\infty.$$
Its closed subspace consisting of holomorphic functions in $\mathcal
D$ is the weighted Bergman space $A^p_\nu(\mathcal D)$. This space
is not trivial i.e. $A^p_\nu(\mathcal D)\neq \{0\}$ only for
$\nu>\frac{n}{r}-1$ (see \cite{DD}). The weighted Bergman projection $P_\nu$ is the orthogonal projection of the Hilbert space $L^2_\nu(\mathcal D)$ onto its closed space $A^2_\nu(\mathcal D).$ It is well known that $P_\nu$ is an integral operator given by $$P_\nu f(z)=\int_{\mathcal D}K_\nu(z, w) f(w)
dV_\nu(w),
$$
where
 $K_\nu(z, w)=
 c_\nu\,\Delta^{-(\nu+\frac{n}{r})}((z-\overline {w})/i)$
is the weighted Bergman kernel, i.e. the reproducing kernel of  $A^2_\nu(\mathcal D)$ (see
\cite{FK}). Here, we use the notation
$dV_\nu(w):=\Delta^{\nu-\frac{n}{r}}(v) du\,dv$, where
$w=u+iv$ is an element of $\mathcal D$. The unweighted case corresponds to $\nu=\frac{n}{r}.$

Given $1\le p,q<\infty$, the first question we consider here is the characterization of those positive measures $\mu$ on $\mathcal D$ such that the embedding $I_\mu:A^p_\nu(\mathcal D)\rightarrow L^q(\mathcal D, d\mu)$ is continuous. That is, there exists a constant $C>0$ such that the following inequality
$$\int_{\mathcal D}|f(z)|^qd\mu(z)\le C\|f\|_{p,\nu}^q$$
holds for any $f\in A^p_\nu(\mathcal D)$, $\nu>\frac{n}{r}-1$. When this holds, we speak of Carleson embedding (or estimate) for $A^p_\nu(\mathcal D)$; we also say that $\mu$ is a $q$-Carleson measure for $A^p_\nu(\mathcal D)$.
\vskip .2cm
Carleson measures play an important role in several problems in mathematical analysis and its applications. Among these problems, we have the questions of the boundedness of Composition operators and Toeplitz operators to name a few.

Carleson estimates without loss ($1\le p\le q<\infty$) for usual (weighted) Bergman spaces of the unit disc or the unit ball have been obtained in \cite{CW, Gu, Hastings, Luecking3}. Estimations with loss ($1\le q<p<\infty$) are essentially due to Luecking \cite{Luecking1}. We also mention the work of O. Constantin \cite{Constantin} in which she considered Bergman spaces with B\'ekoll\'e-Bonami weights.

In this paper, we give a full characterization of Carleson estimates for Bergman spaces in the setting of tube domains over symmetric cones. For this, we will need among others, an atomic decomposition of functions in Bergman spaces with change of weights and in the case of estimations with loss, we shall then make use of the techniques of Luecking \cite{Luecking1}. Unlike the classical case (in the unit disc or the unit ball) where the Bergman projection is bounded on the full range $1<p<\infty$, for tube domains over symmetric cones, the Bergman projection is not bounded for large values of the exponent (see \cite{BBGR} ). Hence our Carleson embeddings with loss are based on the hypothesis that boundedness of the Bergman projection holds, if necessary, with a change of weights.  Next, we use this result to obtain boundedness criteria for the Ces\`aro-type operator defined below.

\vskip .2cm
 In order to  give a definition of  Ces\`aro-type operators adapted to our context, let us go back to the unit disc $\mathbb D$ of $\mathbb C.$  For any $g\in \mathcal {H}(\mathbb D),$ the Ces\`aro-type integral operator $T_g$  was introduced as follows
 $$T_gf(z):=\int_0^z f(w)g'(w)dw,\,\,\,f\in \mathcal {H}(\mathbb D).$$
 It is easy to check that for $f,g\in \mathcal {H}(\mathbb D)$,
 \begin{eqnarray}\label{ces}(T_gf)'(z)=f(z)g'(z).\end{eqnarray}

The operator $T_g$ has been first studied on Bergman spaces of the unit disc by A. Aleman and A. G. Siskakis \cite{AlemanSis}; their work was followed by a heavy literature on this operator and its siblings. We refer to the following and the references therein \cite{Aleman, AlemanConst, Constantin, Siskakis, SisZhao, stevic,StevicUeki, Xiao}.

Recall that the Box operator $\Box=\Delta(\frac{1}{i}
\frac{\partial}{\partial x})$ is the differential operator of
degree $r$ in $\mathbb R^n$ defined by the equality: \Be
\Box\,[e^{i(x|\xi)}]=\Delta(\xi)e^{i(x|\xi)}, \quad
x\in \mathbb R^n,\,\xi\in\O. \label{bbox} \Ee Let us consider $$\mathcal {N}_n := \{F \in
\mathcal H (\mathcal D) : \Box^n F = 0 \}$$ and set
$$\mathcal{H}_n(\mathcal D)=\mathcal{H}(\mathcal D)/ \mathcal N_n.$$
 Taking advantage of the observation (\ref{ces}), we define for $g\in \mathcal {H}(\mathcal D)$ the operator $T_g$
 as follows: for $f\in \mathcal {H}(\mathcal D)$, $T_gf$ is the equivalence class of the solutions of the equation
 $$\Box^{n}F=f\Box^{n}g.$$ Note that the definition of our Ces\`aro-type operator does not depend on the choice of the representative of
the class of the symbol.
We shall use Carleson embeddings to characterize functions $g$ such that $T_g$ extends as a bounded operator from $A_\nu^p(\mathcal D)$ to $\mathcal {B}_\nu^q(\mathcal D)$, where
$\mathcal {B}_\nu^q(\mathcal D)$ is defined as the space of equivalent classes in $\mathcal H_n(\mathcal D)$ of analytic functions $f$ such that
$z\mapsto\Delta^n(\Im m\, z)\Box^nf(z)\in L_\nu^q(\mathcal D)$.

The choice of $\mathcal {B}_\nu^q(\mathcal D)$ as target space is suggested by the choice of the number of derivatives in the definition of $T_g$ and the lack of boundedness
of the Bergman projection.

\vskip .2cm
Our other interest will be for the Schatten class criteria of the Toeplitz operators $T_\mu.$ Here $\mu$ is a positive Borel measure on $\mathcal D.$ For a function $f$ with compact support we define
\begin{equation}\label{defToeplitz}
T_\mu f(z):=\int_{\mathcal D}K_\nu(z,w)f(w)d\mu (w),
\end{equation}
where $K_\nu$ is the weighted Bergman kernel.

Schatten class criteria of the Toeplitz operators have been obtained by several authors
 on bounded domains of  $\mathbb {C}^n$ \cite{Constantin,Luecking2, Zhu2, Zhu3}. More precisely, denote by $\mathbb B^n$ the unit ball of $\mathbb C^n.$ For $\alpha>-1$, write $dv_\alpha(z)=c_\alpha(1-|z|^2)^\alpha dv(z)$ where $c_\alpha$ is such that $v_\alpha(\mathbb {B}^n)=1$, and $dv$ is the Lebesgue measure on $\mathbb B^n$. The Bergman space $A_\alpha^p(\mathbb B^n)$ is defined as above with $\mathcal D$ replaced by $\mathbb B^n$ and $dV_\nu$ by $dv_\alpha$. The corresponding Bergman kernel function is given by $K_\alpha(z,w):=\frac{1}{(1-\langle z,w\rangle)^{n+1+\alpha}}$, where $\langle z,w\rangle=z_1\overline {w_1}+\cdots+z_n\overline {w_n}$, whenever $z=(z_1,\cdots,z_n)$ and $w=(w_1,\cdots,w_n)$.
 \vskip .1cm
 Given a positive measure $\mu$ on $\mathbb B^n$, the Berezin transform $\tilde {\mu}$ of $\mu$ is the function defined on $\mathbb B^n$ by
 $$\tilde {\mu}(z):=\int_{\mathbb B^n}\frac{(1-|z|^2)^{n+1+\alpha}}{|1-\langle z,w\rangle|^{2(n+1+\alpha)}}d\mu(w),\,\,\,z\in \mathbb B^n.$$
 For $z\in \mathbb B^n$ and $\delta>0$ define the average function $$\hat {\mu}_\delta(z):=\frac{\mu(D(z,\delta))}{v_\alpha(D(z,\delta))}$$
 where $D(z,\delta)$ is the Bergman metric ball centered at $z$ with radius $\delta$. The following was proved by K. Zhu in \cite{Zhu2, Zhu3}.
 \btheo\label{theo:zhu} Suppose that $\mu$ is a positive measure on $\mathbb B^n$, $0<p<\infty$, $\delta>0$, and $\alpha>-1$. Then the following assertions are equivalent.
 \begin{itemize}
 \item[(a)] The Toeplitz operator $T_\mu$ belongs to the Schatten class $\mathcal S_p(A_\alpha^2(\mathbb B^n)).$
 \item[(b)] The average function $\hat {\mu}_\delta$ belongs to $L^p(\mathbb B^n, d\lambda)$,\,\,where $d\lambda(z)=\frac{dv(z)}{(1-|z|^2)^{n+1}}$.
 \item[(c)] For any $\delta$-lattice $\{a_k\}$ in the Bergman metric of $\mathbb B^n$, the sequence $\{\hat {\mu}_\delta(a_k)\}$ belongs to $l^p$.

 Moreover, if $p>\frac{n}{n+1+\alpha}$, then the above assertions are equivalent to the following:
 \item[(d)] the Berezin transform $\tilde {\mu}$ belongs to $L^p(\mathbb B^n, d\lambda)$.
 \end{itemize}
 \etheo
Note that the equivalence $(a)\Leftrightarrow (c)$ in the case of the unit disc of $\mathbb C$ is due to Luecking \cite{Luecking2} and that the result in \cite{Zhu2} is in fact stated for Bergman spaces of bounded symmetric domains.

 Our aim here is to formulate and prove these results in the setting of tube domains over symmetric cones. Methods used to study Schatten classes of Toeplitz operators are somehow classical. Ideas used here are analogues of those in \cite{Constantin,Luecking2, S, Choe} and the references therein.

Next, we will consider criteria for Schatten classes membership of the Ces\`aro-type operator above on the weighted Bergman space $A_\nu^2(\mathcal D)$. Our results are in the spirit of the classical ones (see \cite{Constantin, Luecking2, Zhu2, Zhu3} and the references therein). More precisely, we will prove that the range of symbols $g$ for which $T_g$ is in the Schatten class $\mathcal {S}_p$ ($2\le p\le \infty$)
is the Besov space $\mathcal {B}^{p}(\mathcal D)=\mathcal {B}_\nu^p(\mathcal D)$, with $\nu=-n/r$, and for $p=\infty$, $\mathcal {B}^\infty(\mathcal D)=\mathcal {B}_{\nu}^\infty(\mathcal D)$ is the set of all holomorphic functions $f$ on $\mathcal D$ such that $$\sup_{z\in \mathcal D}\Da^n(\im z)|\Box^nf(z)|<\infty.$$ To prove this result, we shall identify our operator with a Toeplitz operator and deduce the Schatten class membership as in classical cases (see \cite{Constantin} for example).

 We remark that criteria for Schatten class membership of the Hankel operators on $A_\nu^2(\mathcal D)$
were obtained by the second author in \cite{S}.

\vskip .2cm
We are mainly motivated by the idea of finding more applications of Besov spaces of this setting as developed in \cite{BBGRS}. We also note that even if many of our ideas (related to the study of Toeplitz and Ces\`aro-type operators in particular), are built from the classical cases, nothing is given for free in this setting as the question of boundedness of the Bergman projection shows. One of the tools we will be using to extend some useful (classical) results in the study of the above operators is the Kor\'anyi lemma (Lemma \ref{lem:Koranyi} in the text).

\vskip .2cm
The paper is organized as follows. In the next section, we present some useful tools needed to prove our results. In section 3, we characterize Carleson measures for weighted Bergman spaces and apply them to the boundedness of our Ces\`aro-type operators. The last section is devoted to Schatten classes criteria of Toeplitz operators and Ces\`aro-type operators.

\vskip .2cm
As usual, given two positive quantities $A$ and $B$, the notation $A\lesssim B$ (resp. $A\gtrsim B$) means that there is an absolute
positive constant $C$ such that $A\le CB$ (resp. $A\ge CB$). When $A\lesssim B$ and $B\lesssim A$, we write $A\simeq B$ and say $A$
and $B$ are equivalent. Finally, all over the text,  $C$, $C_k$, $C_{k,j}$ will
denote positive constants depending only on the displayed
parameters but not necessarily the same at distinct
occurrences.

\section{Preliminary results}

In this section, we give some fundamental facts about symmetric cones, Berezin transform and related results.

\subsection{Symmetric cones and Bergman metric}
Recall that a symmetric cone $\Omega$ induces in $V\equiv \mathbb R^n$
a structure of Euclidean Jordan algebra, in which $\overline
{\Omega}=\{x^2:x\in V\}$. We denote by $e$ the identity element in
$V$. Let $G(\Omega)$ be the group of transformations of $\mathbb R^n$ leaving invariant the cone $\Omega$. The cone $\Omega$ is homogeneous since the
group $G(\Omega)$ acts transitively on $\Omega.$ It is well known that there is a subgroup $H$ of $G(\Omega)$ that acts simply transitively on $\Omega$, that is for $x,y\in\Omega$ there is a unique $h\in H$ such that $y=hx.$ In particular, $\Omega\equiv H\cdot e$.

If we denote by $\mathbb R^n$ the group of translation by vectors in $\mathbb R^n$, then the group $G(\mathcal D)=\mathbb {R}^n\times H$ acts simply transitively on $\mathcal D$.

\vskip .2cm
Next we recall the definition of the Bergman metric. Define a matrix function $\{g_{jk}\}_{1\leq j,k\leq n}$ on $\mathcal D$ by
$$g_{jk}(z)=\frac{\partial^2}{\partial z_j\partial\bar z_k}\log K(z,z)$$
where $K(w,z)=c(n/r)\Delta^{-\frac{2n}{r}}(\frac{w-\overline {z}}{i})$ is the (unweighted) Bergman kernel of $\mathcal D$. The map $z\in \mathcal {D}\mapsto \mathcal G_z$ with
$$\mathcal G_z(u,v)=\sum_{1\leq j,\,k\leq n}g_{jk}(z)u_j{\bar v}_k,\quad u=(u_1,\ldots,u_n),\,\,v=(v_1,\dots,v_n)\in \mathbb C^n$$
defines a Hermitian metric on $\mathbb C^n,$ called the Bergman metric. The Bergman length of a smooth path $\gamma:[0,1]\to \mathcal D$ is given by
$$l(\gamma)=\int_0^1\left(\mathcal G_{\gamma(t)}(\dot\gamma(t),\dot\gamma(t))\right)^\frac{1}{2}dt$$
and the Bergman distance $d(z_1,z_2)$ between two points $z_1, z_2\in \mathcal D$ is defined by
$$d(z_1,z_2)=\inf_\gamma l(\gamma)$$
where the infimum is taken over all smooth paths $\gamma:[0,1]\to T_\Omega$ such that $\gamma(0)=z_1$ and $\gamma(1)=z_2.$

For $\delta>0$, we denote by $$B_\delta(z)=\{w\in \mathcal {D}: d(z,w)<\delta\}$$ the Bergman ball centered at $z$ with radius $\delta.$

We will need the following (see \cite[Theorem 5.4]{BBGNPR}).
\blem\label{lem:covering}
Given $\delta\in (0, 1)$, there exists a sequence $\{z_j\}$ of points  of $\mathcal D$ called
$\delta$-lattice such that, if $B_j=B_\delta(z_j)$ and $B_j'=B_{\frac{\delta}{2}}(z_j)$,
then
\begin{itemize}
\item[(i)] the balls $B_j'$ are pairwise disjoint;
\item[(ii)] the balls $B_j$ cover $\mathcal D$ with finite overlapping, i.e. there is an integer N (depending only on $\mathcal{D}$) such
that each point of $\mathcal D$ belongs to at most N of these balls.
\end{itemize}
\elem

The above balls have the following
properties:$$\int_{B_j}dV_\nu(z)\approx \int_{B'_j}dV_\nu(z)\approx
C_{\delta}\Delta^{\nu+n/r}(\im z_j).$$ We recall that the
measure $d\lambda(z)=\Delta^{-2n/r}(\im z)dV(z)$ is an
invariant measure on $\mathcal D$ under the actions of $G(\mathcal D)=\mathbb {R}^n\times H$.

\vskip .2cm
Let us denote by $l_{\nu}^{p}$, the space of complex sequences $\beta=\{\beta_j\}_{j\in \N}$ such that
$$||\beta||_{l_\nu^p}^{p}=\sum_{j}
|\beta_{j}|^{p}\Delta^{\nu+\frac{n}{r}}(\im z_j)<\infty,$$
where $\{z_j\}_{j\in \N}$ is a $\delta$-lattice.
\blem\label{5.1}\cite[Lemma
2.11.1 ]{G} Suppose $1 \le p < \infty$, $\nu $ and $\mu $
are reals. Then, the dual space $(l_{\nu}^{p})^{*}$ of the space
$l_{\nu}^{p}$ identifies with $l_{\nu + (\mu - \nu)p'}^{p'}$
under the sum pairing $$\langle\eta,\beta\rangle_{\mu}=
\sum_{j}\eta_{j}\overline {\beta_{j}}\Delta^{\mu +
\frac{n}{r}}(\im z_j),$$ where $\eta = \{\eta_{j}\}$ belongs to
$l_{\nu}^{p}$ and $\beta =\{\beta_{j}\}$ belongs to $l_{\nu + (\mu
- \nu)p'}^{p'}$ with $\frac{1}{p} +\frac{1}{p'}=1$.
\elem

We refer to \cite[Theorem 5.6]{BBGNPR} for the following result known as the sampling theorem.
\blem\label{5.3}
Let $\{z_{j}\}_{j\in \N}$ be a $\da$-lattice in $\mathcal D$,
$\da \in (0,1)$ with $z_{j}= x_{j} + iy_{j}$. The following
assertions hold.
\begin{itemize}
\item[(1)] There is a positive
constant $C_{\da}$ such that every $f\in A_{\nu}^{p}(\mathcal D)$ satisfies
$$ ||\{f(z_{j})\}||_{l_\nu^p} \le C_{\da} ||f||_{p,\nu}. $$
\item[(2)] Conversely, if $\da$ is small enough, there is a
positive constant $C_{\da}$ such that every $f\in A_{\nu}^{p}(\mathcal D)$
satisfies $$||f||_{p,\nu} \le
C_{\da}||\{f(z_{j})\}||_{l_\nu^p}.$$
\end{itemize}
\elem

We will need the following consequence of the mean value theorem (see \cite{BBGNPR})
\blem\label{lem:meanvalue}
There exists a constant $C>0$ such that for any $f\in \mathcal {H}(\mathcal D)$ and $\delta\in (0,1]$, the following holds
\begin{equation}\label{eq:meanvalue}
|f(z)|^p\le C\delta^{-n}\int_{B_\delta(z)}|f(\zeta)|^p\frac{dV(\zeta)}{\Delta^{2n/r}(\im \zeta)}.
\end{equation}
\elem

We finish this subsection by recalling that there is a constant $C>0$ such that for any $f\in A_\nu^p(\mathcal D)$ the following
pointwise estimate holds:
\begin{equation}\label{eq:pointwiseestimberg}|f(z)|\le C\Delta^{-\frac{1}{p}(\nu+\frac{n}{r})}(\im z)\|f\|_{A_\nu^p},\,\,\,\textrm{for all}\,\,\, z\in \mathcal {D}
\end{equation}
(see \cite[Proposition 3.5]{BBGNPR}).
\subsection{Averaging functions and Berezin transform}
Let $\mu$ be a positive measure on $\mathcal D$. For $\nu>\frac{n}{r}-1$ and $w\in \mathcal D$, for simplicity, let us consider the normalized reproducing kernel
\begin{equation}\label{eq:normalrepkern} k_\nu(\cdot,w)=\frac{K_\nu(\cdot,w)}{\|K_\nu(\cdot,w)\|_{2,\nu}}=\Da^{-\nu-\frac nr}\left(\frac{\cdot-\bar w}{i}\right)\Da^{\frac{1}{2}(\nu+\frac nr)}(\im w).
\end{equation}

For $z\in \mathcal D$, we define
$$\tilde {\mu}(z):=\int_{\mathcal D}|k_\nu(z,w)|^2d\mu(w).$$
The function $\tilde {\mu}$ is the Berezin transform of the measure $\mu$. When $d\mu(z)=f(z)dV_\nu(z)$, we write $\tilde {\mu}=\tilde {f}$ and speak of the Berezin transform of the function $f$. 
The following is given in \cite[Theorem 1.1]{sehba}.
\blem\label{lem:Berenzinbounded}
Let $\nu>\frac{n}{r}-1$, $\alpha\in \mathbb R$, and $1<p<\infty$. Then the operator $f\mapsto \tilde {f}$ (associated to $\nu$) is bounded on $L_\alpha^p(\mathcal D)$ if and only if 
$$\nu p-\alpha>\left(\frac{n}{r}-1\right)\max\{1,p-1\}$$
and $$(\nu+\frac{n}{r})p+\alpha>\left(\frac{n}{r}-1\right)\max\{1,p-1\}.$$
\elem

For $z\in \mathcal D$ and $\delta\in (0,1)$, we define the average of the positive measure $\mu$ at $z$ by $$\hat {\mu}_\delta(z)=\frac{\mu(B_\delta(z))}{V_\nu(B_\delta(z))}.$$
The function $\tilde {\mu}$ and $\hat {\mu}$ are very useful in the characterization of Schatten classes for Toeplitz operators (see \cite{Zhu2,Zhu3} and the references therein).

We have the following result known as Kor\'anyi's lemma.
\blem\label{lem:Koranyi}\cite[Theorem 1.1]{BIN}
For every $\delta>0,$ there is a constant $C_{\delta}>0$ such that
$$\left|\frac{K(\zeta,z)}{K(\zeta,w)}-1\right|\leq C_{\delta}d(z,w)$$
for all $\zeta,z,w\in \mathcal D,$ with $d(z,w)\leq \delta.$
\elem
The following corollary is a straightforward consequence of \cite[Corollary 3.4]{BBGNPR} and Lemma \ref{lem:Koranyi}.
\bcor\label{kor}
Let $\nu>\frac nr-1,\,\,\da>0$ and $z,w\in \mathcal D.$  
There is a positive constant $C_\da$ such that for all $z\in B_\da(w),$
$$V_\nu(B_\da(w))|k_\nu(z,w)|^2\leq C_\da.$$

If $\da$ is sufficiently small, then there is $C>0$ such that for all $z\in B_\da(w),$
$$V_\nu(B_\da(w))|k_\nu(z,w)|^2\geq (1-C\da).$$
\ecor
The following two results are direct consequences of the above corollary.
\blem\label{lem:Berezinsupaverage}
Let $\delta\in (0,1)$. Then there exists a constant $C=C_\delta>0$ such that
$$\hat {\mu}_\delta(z)\le C_\delta\tilde {\mu}(z),\,\,\,\textrm{for any}\,\,\,z\in \mathcal D.$$
\elem
\begin{proof}
Let $z\in \mathcal D$ be given. Using Corollary \ref{kor} we obtain
$$\int_{B_\delta(z)}|k_\nu(z,\xi)|^2d\mu(\xi)\ge C_\delta \frac{\mu(B_\delta(z))}{V_\nu(B_\delta(z))}.$$
Thus $\tilde {\mu}(z)\ge C_\delta\hat {\mu}_\delta(z)$ for any $z\in \mathcal D$.
\end{proof}

\bcor\label{lem:equivVolume}
Let $\nu>\frac{n}{r}-1$, $\delta,t\in (0,1)$. Let $0<\delta_1,\delta_2,\delta_3<\delta$ with $t<\delta_1/ \delta_2<t^{-1}$. Then there exists a constant $C=C(\delta,t)>0$ such that for any $z,w\in \mathcal D$ such that $w\in B_{\delta_3}(z)$, $$\frac{1}{C}<\frac{V_\nu(B_{\delta_1}(z))}{V_\nu(B_{\delta_2}(w))}<C.$$
\ecor
The following lemma will allow flexibility on the choice of the radius of the ball for some questions along the text. It is adapted from \cite{Choe}.
\blem\label{lem:variationoflattice}
Let $1\le p\le \infty$, $\nu\in \mathbb R$, and $\delta,\beta\in (0,1)$. Assume that $\mu$ is a positive Borel measure on $\mathcal D$. Then
the following assertions are equivalent.
\begin{itemize}
\item[(i)] The function  $\mathcal {D}\ni z\mapsto \frac{\mu(B_\delta(z))}{\Delta^{\nu+\frac{n}{r}}(\im z)}$ belongs to $L_\nu^p(\mathcal D)$.
\item[(ii)] The function  $\mathcal {D}\ni z\mapsto \frac{\mu(B_\beta(z))}{\Delta^{\nu+\frac{n}{r}}(\im z)}$ belongs to $L_\nu^p(\mathcal D)$.
\end{itemize}
\elem
\Proof
We only need to prove that $(i)\Rightarrow (ii)$. We may suppose that $\delta<\beta$ as otherwise the implication is obvious. Choose a finite set $\{z_1,\cdots,z_m\}\subset B_\beta(ie)$ which is maximal with respect to condition $B_{\delta /2}(z_j)\cap B_{\delta /2}(z_k)=\emptyset$ for $j\neq k$.
Then $B_{\beta}(ie)\subset \cup_{j=1}^mB_{\delta }(z_j)$ by the maximality property of $\{z_1,\cdots,z_m\}$ and consequently, for any $z\in \mathcal D$ and for some appropriate $g\in G(\mathcal D)$,
$$B_\beta(z)=g\left(B_\beta(ie)\right)\subset \cup_{j=1}^mg\left(B_\delta(z_j)\right)=\cup_{j=1}^mB_\delta(gz_j).$$
We recall that as $gz_j\in B_\beta(z)$ and $\delta<\beta$, we have $B_\delta(gz_j)\subset B_{2\beta}(z).$ Therefore, $$V_\nu(B_\beta(z))\leq \sum_{j=1}^mV_\nu(B_\delta(gz_j))\leq mV_\nu(B_{2\beta}(z))\simeq V_\nu(B_\beta(z))$$ for all $z\in \mathcal D$. Hence, $V_\nu(B_\delta(gz_j))\leq C V_\nu(B_{\beta}(z)).$ It follows that
\Beas
\frac{\mu(B_\beta(z))}{V_\nu\left(B_\beta(z)\right)} &\le& \sum_{j=1}^m\frac{\mu(B_\delta(gz_j))}{V_\nu\left(B_\beta(z)\right)}\\ &=&
\sum_{j=1}^m\frac{\mu(B_\delta(gz_j))}{V_\nu\left(B_\delta(gz_j)\right)}\times \frac{V_\nu\left(B_\delta(gz_j)\right)}{V_\nu\left(B_\beta(z)\right)}\\ &\leq& C \sum_{j=1}^m\frac{\mu(B_\delta(gz_j))}{V_\nu\left(B_\delta(gz_j)\right)}.
\Eeas
This proves that $(i)\Rightarrow (ii)$ for $p=\infty$.

For $p$ finite, the implication follows from the following inequality:
$$\left(\frac{\mu(B_\beta(z))}{V_\nu\left(B_\beta(z)\right)}\right)^p\lesssim m^{p-1}\sum_{j=1}^m\left(\frac{\mu(B_\delta(gz_j))}{V_\nu\left(B_\delta(gz_j)\right)}\right)^p.$$
To finish the proof, observe that for each $z=x+iy\in \mathcal D,$ since $\Omega=He$, we write $y=g(y)e$ with $g(y)\in G(\Omega)$.
The affine linear transformation $$g_z:w\mapsto g_z(w)=g(\im z)w+\Re e\, z$$
is in $G(\mathcal D)$ (see \cite[Proposition X.5.4]{FK}) and satisfies $g_z(ie)=z$.

The conclusion follows now from the following proposition.
\bprop
Let $w\in \mathcal D$ be fixed, $f$ be a measurable function in $\mathcal D$, and $\nu\in \mathbb R$. Then
\begin{equation}
\int_{\mathcal D}|f(g_zw)|dV_\nu(z)=\Delta^{-\nu-\frac{n}{r}}(\im w)\int_{\mathcal D}|f(\zeta)|dV_\nu(\zeta).
\end{equation}
\eprop
\Proof
We use the change of variables $z=h_{zw}\zeta$ where $h_{zw}=g_zg_w^{-1}g_z^{-1}\in G(\mathcal D)$. We have ${\rm Det} h_{zw}=[{\rm Det} g(\im w)]^{-1}$ where $g(\im w)\in G(\O).$  Hence from \cite[Proposition III.4.3]{FK},
\Beas
dV(z) &=& |{\rm Det} h_{zw}|^2dV(\zeta)= \left({\rm Det} g(\im w)\right)^{-2}dV(\zeta)= \Delta^{-2n/r}(\im w)dV(\zeta).
\Eeas
and since $\im z=g(\im z)g(\im w)^{-1}g(\im z)^{-1}\im \zeta,$ we have
\Beas
\Delta^{\nu-\frac nr}(\im z) = \left({\rm Det} g(\im w)\right)^{-\frac{r}{n}(\nu-\frac nr)}\Delta^{\nu-\frac nr}(\im \zeta)= \Delta^{-\nu+\frac{n}{r}}(\im w)\Delta^{\nu-\frac nr}(\im \zeta).
\Eeas
It follows that
\Beas
dV_\nu(z) &=& \Delta^{\nu-\frac nr}(\im z)dV(z) = \Delta^{-\nu-\frac nr}(\im w)dV_\nu(\zeta).
\Eeas
Thus
$$\int_{\mathcal D}|f(g_zw)|dV_\nu(z)=\Delta^{-\nu-\frac{n}{r}}(\im w)\int_{\mathcal D}|f(\zeta)|dV_\nu(\zeta).$$
\ProofEnd
\ProofEnd
We will be sometimes making use of the above lemma in the next sections without any further reference to it.

We also observe the following.
\blem\label{lem:AveragesupBerezin}
Given $\delta\in (0,1)$, there exists $C_\delta>0$ such that
$$\tilde {\mu}(z)\le C_\delta\Tilde {\hat {\mu}}_\delta(z),\,\,\,\textrm{for any}\,\,\,z\in \mathcal D.$$
\elem
\begin{proof}
First observe that by Lemma \ref{lem:meanvalue}, we can write
\Beas
|f(z)|^p\leq\frac{1}{V_\nu(B_\da(z))}\int_{B_\da(z)}|f(w)|^pdV_\nu(w).
\Eeas
It follows with the help of Fubini's lemma that
\Beas
\tilde {\mu}(z) &=& \int_{\mathcal D}|k_\nu(z,\xi)|^2d\mu(\xi)\\ &\lesssim& \int_{\mathcal D} \frac{1}{V_\nu(B_\da(\xi))}\int_{B_\da(\xi)}|k_\nu(z,w)|^2dV_\nu(w)d\mu(\xi)\\ &=& \int_{\mathcal D} \frac{1}{V_\nu(B_\da(\xi))}\int_{B_\da(\xi)}|k_\nu(z,w)|^2\chi_{B_\delta(\xi)}(w)dV_\nu(w)d\mu(\xi)\\ &=& \int_{\mathcal D} \frac{1}{V_\nu(B_\da(\xi))}\int_{B_\da(\xi)}|k_\nu(z,w)|^2\chi_{B_\delta(w)}(\xi)dV_\nu(w)d\mu(\xi)\\ &\backsimeq& \int_{\mathcal D}\frac{\mu(B_\delta(w))}{V_\nu(B_\delta(w))}|k_\nu(z,w)|^2dV_\nu(w)\\ &=& \Tilde {\hat {\mu}}_\delta(z).
\Eeas
We have used that $\chi_{B_\delta(\xi)}(w)=\chi_{B_\delta(w)}(\xi)$ and that as $w\in B_\delta(\xi)$, $V_\nu(B_\delta(\xi))\backsimeq V_\nu(B_\delta(w))$.
\end{proof}
We now prove the following useful result.
\blem\label{lem:integraldiscretizationAverBer}
Let $1\le p\le \infty$, $\nu>\frac{n}{r}-1$, $\alpha\in \mathbb R$, $\beta, \delta\in (0,1)$. Let $\{z_j\}_{j\in \mathbb N}$ be a $\delta$-lattice in $\mathcal D$, and let $\hat {\mu}_\beta$ and $\tilde {\mu}$ be in this order, the average function and the Berezin transform associated to the weight $\nu$ . Then the following assertions are equivalent.
\begin{itemize}
\item[(i)] $\hat {\mu}_\beta\in L_\alpha^p(\mathcal D)$.
\item[(ii)] $\{\hat {\mu}_\delta(z_j)\}_{j\in \mathbb N}\in l_\alpha^p$.

Moreover, if \Be\label{eq:Berezincond1}\nu p-\alpha>\left(\frac{n}{r}-1\right)\max\{1,p-1\}\Ee
and \Be\label{eq:Berezincond2}(\nu+\frac{n}{r})p+\alpha>\left(\frac{n}{r}-1\right)\max\{1,p-1\},\Ee 
then the above assertions are equivalent to
\item[(iii)] $\tilde {\mu}\in L_\alpha^p(\mathcal D)$.
\end{itemize}
\elem
\begin{proof}
The proof is essentially adapted from the unit ball version in \cite{Choe}. We prove that $(i)\Leftrightarrow (ii)$ and $(i)\Leftrightarrow (iii)$. We start by the latter.

Suppose that $\nu>\frac{n}{r}-1$, and (\ref{eq:Berezincond1}) and (\ref{eq:Berezincond2}) hold. That $(iii)\Rightarrow (i)$ is direct from Lemma \ref{lem:Berezinsupaverage}. That $(i)\Rightarrow (iii)$ is direct from Lemma \ref{lem:AveragesupBerezin}
and Lemma \ref{lem:Berenzinbounded}.

$(i)\Rightarrow (ii)$: By Lemma \ref{lem:variationoflattice}, it is enough to prove the estimate
\begin{equation}\label{eq:hatestim}
\|\{\hat {\mu}_\delta(z_j)\}\|_{l_\alpha^p}\lesssim \|\hat {\mu}_{\delta+\gamma}\|_{L_\alpha^p}
\end{equation}
for some $\gamma<\min\{\frac{\delta}{2},1-\delta\}$.

Recall that the balls $B_{\delta/2}(z_j)$ are pairwise disjoint; consequently, the balls $B_\gamma(z_j)$ are also pairwise disjoint. Also, note that for $w\in B_{\gamma}(z_j)$, $B_{\delta}(z_j)\subset B_{\delta+\gamma}(w)$. Hence
\Bea\label{delta1}
\hat {\mu}_{\delta+\gamma}(w) &\ge& \frac{\mu(B_\delta (z_j))}{V_\nu(B_{\delta+\gamma}(w))}\nonumber\\ &=& \frac{V_\nu(B_\delta(z_j))}{V_\nu(B_{\delta+\gamma}(w))}\hat {\mu}_\delta(z_j)\nonumber\\ &\backsimeq& \hat {\mu}_\delta(z_j).
\Eea
It follows that (\ref{eq:hatestim}) holds for $p=\infty$. For $1\le p<\infty$, we obtain using (\ref{delta1})
\Beas
\|\{\hat {\mu}(z_j)\}\|_{l_\alpha^p}^p &\backsimeq& \sum_j\left(\hat {\mu}_\delta(z_j)\right)^pV_\alpha(B_{\delta/2}(z_j))\\ &\lesssim& \sum_{j}\int_{B_{\delta/2}(z_j)}\hat {\mu}_{\delta+\gamma}(w)^pdV_\alpha(w)\\ &\le& \int_{\mathcal D}\hat {\mu}_{\delta+\gamma}(w)^pdV_\alpha(w).
\Eeas
\vskip .2cm
$(ii)\Rightarrow (i)$: Following Lemma \ref{lem:variationoflattice}, we may assume that $\beta<1-\delta$. Given $w\in \mathcal D$, consider the set $$N(w):=\{j:B_\delta(z_j)\cap B_\beta(w)\neq \emptyset\}.$$

The set $N(w)\neq\emptyset$ because $\{z_j\}_{j\in \mathbb N}$ is a $\delta$-lattice on $\mathcal D$ and so the balls $B_\delta(z_j)$ cover $\mathcal D.$ Thus, $$B_\beta(w)\subset \bigcup_{j\in N(w)}B_\delta(z_j).$$
It follows that
\begin{equation}\label{eq:hatbetadeltaestim}
\hat {\mu}_\beta(w) \le \sum_{j\in N(w)}\frac{V_\nu(B_\delta(z_j))}{V_\nu(B_\beta(w))}\hat {\mu}_\delta(z_j) \backsimeq \sum_{j\in N(w)}\hat {\mu}_\delta(z_j).
\end{equation}
Next we observe that if $|N(w)|$ is the cardinality of the set $N(w),$ then
we have $$N_{\beta\,\delta}:=\sup_{w\in \mathcal D}|N(w)|<\infty.$$

To see this, fix $\gamma<\min\{\frac{\delta}{2}, 1-\beta-\delta\}$. Recall that the balls $B_{\gamma}(z_j)$ are pairwise disjoint and observe that for $j\in N(w)$, they are all contained in $B_{\beta+\delta+\gamma}(w)$. It follows that
\Bea\label{delta2} \bigcup_{j\in N(w)}B_{\gamma}(z_j)\subset B_{\beta+\delta+\gamma}(w).\Eea
Let $j\in N(w),$ then $w\in B_{\gamma+\beta}(z_j).$ Thus
\Bea\label{delta3} V_\nu(B_{\beta+\delta+\gamma}(w))\simeq V_\nu(B_{\gamma+\beta}(z_j))\simeq V_\nu(B_{\gamma}(z_j)).\Eea
Therefore, from (\ref{delta2}) and (\ref{delta3}), we get
$$|N(w)|V_\nu(B_{\beta+\delta+\gamma}(w))\backsimeq \sum_{j\in N(w)}V_\nu(B_{\gamma}(z_j))\le V_\nu(B_{\beta+\delta+\gamma}(w)).$$

Thus (\ref{eq:hatbetadeltaestim}) provides that $$\|\hat {\mu}_\beta\|_{L^\infty}\lesssim \|\{\hat {\mu}_\delta(z_j)\}\|_{l^\infty}.$$
For $1\le p<\infty$, we have from (\ref{eq:hatbetadeltaestim}) that
$$\left(\hat {\mu}_\beta(w)\right)^p\lesssim N_{\beta\,\delta}^{p-1}\sum_{j\in N(w)}\left(\hat {\mu}_\delta(z_j)\right)^p.$$
Integrating both sides of the latter inequality with respect to $dV_\alpha(w)$ and applying Fubini's lemma, we obtain
\Beas
\int_{\mathcal D}\left(\hat {\mu}_\beta(w)\right)^pdV_\alpha(w) &\lesssim& \int_{\mathcal D}\sum_{j\in N(w)}\left(\hat {\mu}_\delta(z_j)\right)^pdV_\alpha(w)\\ &=& \sum_{j\in \mathbb N}\left(\hat {\mu}_\delta(z_j)\right)^pV_\alpha(Q_j)
\Eeas
where $$Q_j=\{w\in \mathcal D:B_\delta(z_j)\cap B_\beta(w)\neq \emptyset\}.$$
But clearly, $Q_j\subset B_{\beta+\delta}(z_j)$; hence $$V_\alpha(Q_j)\le V_\alpha(B_{\beta+\delta}(z_j))\backsimeq \Delta^{\alpha+\frac{n}{r}}(\im z_j).$$
We conclude that
$$\int_{\mathcal D}\left(\hat {\mu}_\beta(w)\right)^pdV_\alpha(w)\lesssim \sum_{j\in \mathbb N}\left(\hat {\mu}_\delta(z_j)\right)^p\Delta^{\alpha+\frac{n}{r}}(\im z_j)<\infty.$$ That is $(i)$ holds and the proof is complete.
\end{proof}

\section{Carleson measures for weighted Bergman spaces and boundedness of Ces\`aro-type operators}
In this section, we study Carleson embeddings for weighted Bergman spaces and as application, we give conditions
on the symbols $g$ such that the Ces\`aro-type operators $T_g$ are bounded from $A_\nu^p(\mathcal D)$ to
$\mathcal {B}_\nu^q(\mathcal D)$, $1\le p,q<\infty$. For $1\le p\le q<\infty$, $q$-Carleson measures for usual (weighted) Bergman spaces of the unit disc or the unit ball have been obtained in \cite{CW, Gu, Hastings, Luecking3}. Estimations with loss are essentially due to Luecking \cite{Luecking1}.  In \cite{Constantin}, O. Constantin dealed with the case of Bergman spaces with B\'ekoll\'e-Bonami weights in the unit disc. We first provide an atomic decomposition of functions in Bergman spaces with change of weights.

\subsection{Atomic Decomposition of weighted Bergman spaces}
We give in this subsection, an atomic decomposition  related to a
$\delta$-lattice in $\mathcal {D}$ for
functions in $A_{\nu}^{p}(\mathcal D)$ for
any $1< p < \infty$.
%

For $1 <p < \infty$, we
consider for all $\nu >\frac{n}{r}-1$ the integral pairing
$$\langle f,g\rangle_{\nu}=\int_{\mathcal D}f(z)\overline {g(z)}\Delta^{\nu
-\frac{n}{r}}(\im z)dV(z)$$
where $f\in L_\nu^p(\mathcal {D})$ and $g\in L_\nu^{p'}(\mathcal {D}).$
It is well known that if $P_{\nu}$ is
bounded on $L_{\nu}^{p}(\mathcal D),$ then the dual space of $A_{\nu}^{p}(\mathcal D)$
identifies with $A_{\nu}^{p'}(\mathcal D)$ under the form $\langle \cdot ,\cdot\rangle_{\nu}$.
\blem\label{5.2} Suppose $\mu,\,\,\nu >\frac{n}{r}-1$, $\nu + (\mu
- \nu)p'
> \frac{n}{r}-1$, $1 < p <
\infty$. Then if $P_{\mu}$ is bounded on $L_{\nu}^{p}(\mathcal D)$, then
the following statements hold
\begin{itemize}
\item[(a)] $P_{\mu}$
is the identity on $A_{\nu + (\mu - \nu)p'}^{p'}$, in
particular $P_{\mu}\left(L_{\nu + (\mu - \nu)p'}^{p'}(\mathcal D)\right)= A_{\nu +
(\mu - \nu)p'}^{p'}$.
\item[(b)] $A_{\nu + (\mu -
\nu)p'}^{p'}(\mathcal D)$ is the dual space of the Bergman space
$A_{\nu}^{p}(\mathcal D)$ relatively to the integral pairing $\langle \cdot,\cdot\rangle_{\mu}$.
\end{itemize}
\elem
\Proof Let us recall that if $P_{\mu}$ is bounded on $L_{\nu}^{p}(\mathcal D)$,
then $P_{\mu}$ is the identity on $A_{\nu}^{p}(\mathcal D)$ and moreover,
$P_{\mu}\left(L_{\nu}^{p}\right)= A_{\nu}^{p}$ (see \cite[Lemma 5.1]{BBPR} ).

\vskip .2cm
$(a)$ Clearly, every element of $L_{\nu + (\mu - \nu)p'}^{p'}(\mathcal D)$ defines
an element of the dual space $\left(L_{\nu}^{p}(\mathcal D)\right)^{*}$ of the Lebesgue
space $L_{\nu}^{p}(\mathcal D)$ with respect to the integral pairing
 $\langle \cdot,\cdot\rangle_{\mu}$ and, the operator $S(f)(x+iy)=\Da^{\mu -
\nu}(y)f(x+iy)$ is an isometric isomorphism from $L_{\nu + (\mu - \nu)p'}^{p'}(\mathcal D)$
onto $L_{\nu}^{p'}(\mathcal D)$ such that for every $f\in
L_{\nu}^{p}(\mathcal D)$ and $g\in L_{\nu + (\mu - \nu)p'}^{p'}(\mathcal D)$, $$\langle f,g\rangle_{\mu} =
\langle f,Sg\rangle_{\nu}.$$
 We conclude that the dual space of the Lebesgue
space $ L_{\nu}^{p}(\mathcal D)$ identifies with $L_{\nu + (\mu -
\nu)p'}^{p'}(\mathcal D)$ under the integral pairing $\langle \cdot,\cdot\rangle_{\mu}$.

 Since $P_{\mu}$  is self-adjoint relatively to the form
 $\langle \cdot,\cdot\rangle_{\mu}$, $P_{\mu}$ is also bounded on $L_{\nu + (\mu - \nu)p'}^{p'}$ and
the results then follow from what we have said at the beginning of
the proof.

$(b)$ Every element of $A_{\nu + (\mu - \nu)p'}^{p'}(\mathcal D)$ belongs
to the dual space of the Bergman space $A_{\nu}^{p}(\mathcal D)$ with
respect to the duality pairing  $\langle \cdot,\cdot\rangle_{\mu}$ and since
$P_{\mu}\left(L_{\nu + (\mu - \nu)p'}^{p'}(\mathcal D)\right)= A_{\nu + (\mu -
\nu)p'}^{p'}(\mathcal D),$  we have that every element of the dual space
$\left(A_{\nu}^{p}(\mathcal D)\right)^{*}$ of the Bergman space $A_{\nu}^{p}(\mathcal D)$
identifies with an element of $A_{\nu + (\mu - \nu)p'}^{p'}(\mathcal D)$.
This completes the proof of the lemma.
\ProofEnd

We have the following atomic decomposition with change of weight.
\btheo\label{theo:atomdecompo} Let $1 < p <\infty$ and let $\mu,\,\,\nu >\frac{n}{r}-1$ satisfying $\nu +
(\mu - \nu)p' > \frac{n}{r}-1.$ Assume that
the operator $P_{\mu}$ is bounded on $L_{\nu}^{p}(\mathcal D)$ and let $\{z_{j}\}_{j\in \N}$
be a $\da$-lattice in $\mathcal D$. Then the following assertions hold.
\begin{itemize}
\item[(i)] For every complex sequence $\{\la_{j}\}_{j\in \N}$ in
$l_{\nu}^{p}$, the series  $\sum_{j} {\lambda_{j}K_{\mu}(z,
z_{j})\Delta^{\mu + \frac{n}{r}}(y_{j})}$ is convergent in
$A_{\nu}^{p}(\mathcal D)$. Moreover, its sum $f$ satisfies the inequality
$$||f||_{p,\nu}\leq C_{\delta}||\{\la_{j}\}||_{l_{\nu}^{p}},$$ where
$C_{\delta}$ is a positive constant.
\item[(ii)] For $\da$ small
enough, every function $f\in A_{\nu}^{p}(\mathcal D)$ may be written as
$$ f(z) = \sum_{j} {\lambda_{j}K_{\mu}(z,
z_{j})\Delta^{\mu + \frac{n}{r}}(y_{j})}$$ with
\begin{equation}\label{eq:reverseineqatomdecomp}
||\{\la_{j}\}||_{l_{\nu}^{p}}\le
C_{\da}||f||_{p,\nu}\Ee where $C_{\da}$ is a positive
constant.
\end{itemize}
\etheo
\Proof The proof follows the same steps as the one for the case $\mu=\nu$ (see \cite{BBGNPR} and the references therein), we give it here as our result is more general and the change of weight matters. We already know from Lemma \ref{5.2} that the dual space
 $\left(A_{\nu}^{p}(\mathcal D)\right)^{*}$ of the Bergman space $A_{\nu}^{p}(\mathcal D)$ identifies with
 $A_{\nu + (\mu - \nu)p'}^{p'}(\mathcal D)$ under the integral pairing
 $$\langle f,g\rangle_{\mu}=\int_{\mathcal D}f(z)\overline {g(z)}\Da^{\mu -\frac{n}{r}}(\im z)dV(z).$$
 We also have from Lemma \ref{5.1}
that the dual space $(l_{\nu}^{p})^{*}$ of $l_{\nu}^{p}$
identifies with $l_{\nu + (\mu - \nu)p'}^{p'}$ under the sum
pairing $$\langle \beta,\alpha\rangle_{\mu}= \sum_{j}\beta_{j}\overline
{\ba_{j}}\Da^{\mu + \frac{n}{r}}(y_{j}).$$

$(i)$ From part $(1)$ of the sampling theorem (Lemma
\ref{5.3}), we deduce that the linear operator
$${R:\,\,\,A_{\nu}^{p}(\mathcal D)\rightarrow l_{\nu}^{p}}$$ $${f\mapsto \{f(z_{j})\}}$$
 is bounded. Hence its  adjoint under
the form $\langle \cdot,\cdot\rangle_{\mu}$\,\,$$R^{*}:\,\,\,l_{\nu + (\mu -
\nu)p'}^{p'}\rightarrow A_{\nu + (\mu - \nu)p'}^{p'}$$ is
also bounded. To conclude, it suffices to show that
$$R^{*}(\{\la_{j}\})= \sum_{j} {\lambda_{j}K_{\mu}(\cdot,
z_{j})\Delta^{\mu + \frac{n}{r}}(y_{j})}.$$
For $f\in
A_{\nu}^{p}(\mathcal D)$ and $\{\la_{j}\}\in l_{\nu + (\mu -
\nu)p'}^{p'}$, we have 
\Beas
 \langle Rf,\{\la_{j}\}\rangle_{\mu}&=&
\sum_{j} f(z_{j})\overline {\la_{j}}\Da^{\mu +
\frac{n}{r}}(y_{j})\\ & =& \sum_{j} (P_{\mu}f(z_{j}))\overline
{\la_{j}}\Da^{\mu + \frac{n}{r}}(y_{j})\\ &=& \sum_{j}
\int_{\mathcal D}K_{\mu}(z_{j},w)f(w)dV_\mu(w)\overline {\la_{j}}\Da^{\mu +
\frac{n}{r}}(y_{j})\\ &=& \int_{\mathcal D}f(w)\overline{\left(\sum_{j}
{\la_{j}}K_{\mu}(w,z_{j})\Da^{\mu +
\frac{n}{r}}(y_{j})\right)}dV_\mu(w)\\ &=&
\langle f,R^{*}(\{\la_{j}\})\rangle_{\mu}\Eeas and the conclusion then
follows.

\vskip .2cm
 $(ii)$ From part $(2)$ of Lemma \ref{5.3}, we have that for
$\da$ small enough,
$$||f||_{p',\nu + (\mu - \nu)p'} \le C_{\da}||\{f(z_{j})\}||_{l_{\nu + (\mu - \nu)p'}^{p'}}. $$
This implies that $R^{*}:\,\,\,l_{\nu + (\mu -
\nu)p'}^{p'}\rightarrow A_{\nu + (\mu - \nu)p'}^{p'}(\mathcal D)$ is
onto. Moreover, if we denote by $\mathcal {N}$ the subspace of
$l_{\nu + (\mu - \nu)p'}^{p'}$ consisting of all sequences
$\{\la_{j}\}_{j\in \N}$ such that the sum $\sum_{j}
{\lambda_{j}B_{\mu}(z, z_{j})\Delta^{\mu +
\frac{n}{r}}(y_{j})}$ is identically zero, then the linear map
$$l_{\nu + (\mu - \nu)p'}^{p'}/\mathcal {N} \rightarrow A_{\nu + (\mu - \nu)p'}^{p'}(\mathcal D)$$
is a bounded isomorphism. The continuity of its inverse which
follows from the Hahn-Banach theorem gives the estimate (\ref{eq:reverseineqatomdecomp}). This
completes the proof of the theorem.
\ProofEnd
\begin{Remark}
Note that for $1<p\le 2$ and $\nu>\frac{n}{r}-1$, an atomic decomposition as above always holds for $A_\nu^p(\mathcal D)$ as for example
taking $\mu=\nu+\frac{n}{r}$, \cite[Theorem 1.1]{sehba} gives that $P_\mu$ is bounded on $L_\nu^p(\mathcal D)$.
\end{Remark}

\subsection{Carleson embeddings for weighted Bergman spaces}
We start this subsection with the following observation.
\bprop\label{prop:Carlembedlattice}
Let $1\le p<\infty$, $\nu>\frac{n}{r}-1$. If
$$\sup_{j}\frac{\mu(B_\delta(\zeta_j))}{\Delta^{\nu+\frac{n}{r}}(\im \zeta_j)}\le C<\infty,$$
for some $\delta\in (0,1)$ where $\{\zeta_j\}$ is a
$\delta$-lattice, then
$$\int_{\mathcal D}|f(z)|^pd\mu (z)\le C\|f\|_{p,\nu}^p,\,\,\,f\in A_\nu^p(\mathcal D).$$
\eprop
\Proof
Clearly, we have using Lemma \ref{lem:covering}
\Beas
\int_{\mathcal D}|f(z)|^pd\mu (z) &\le& C_p\sum_j\int_{B_j}|f(z)|^pd\mu (z)\\ &\leq& C_p \sum_j \left(\sup_{z\in B_j}|f(z)|^p\right)\mu(B_j).
\Eeas
Let $z\in B_j.$  By Lemma \ref{lem:integraldiscretizationAverBer}, we can suppose that $\delta<\frac{1}{2}$. Since $B_\da(z)\subset B_{2\da}(z_j)$ and $V_\nu(B_\da(z))\simeq V_\nu(B_\da(z_j))$ (cf. \cite[Theorem 2.38]{BBGNPR}), by Lemma \ref{lem:meanvalue}, we write
\Beas
|f(z)|^p\leq\frac{1}{V_\nu(B_\da(z))}\int_{B_\da(z)}|f(w)|^pdV_\nu(w)\leq \frac{C}{V_\nu(B_\da(z_j))}\int_{B_{2\da}(z_j)}|f(w)|^pdV_\nu(w).
\Eeas
Hence, $$\sup_{z\in B_j}|f(z)|^p\leq \frac{1}{V_\nu(B_\da(z_j))}\int_{B_{2\da}(z_j)}|f(w)|^pdV_\nu(w)$$
and
\Beas
\int_{\mathcal D}|f(z)|^pd\mu (z) &\le& C_pC\sum_j\frac{\mu(B_j)}{V_\nu(B_\da(z_j))}\int_{B_{2\da}(z_j)}|f(w)|^pdV_\nu(w)\\
&\leq&C'_p\sum_j\int_{B_{2\da}(z_j)}|f(w)|^pdV_\nu(w)\leq C'_pM\int_{\mathcal D}|f(w)|^pdV_\nu(w)
\Eeas
since each point of ${\mathcal D}$ belongs to at most $M$ balls $B_{2\da}(z_j).$
\ProofEnd
Our next result on the Carleson embeddings for weighted Bergman spaces is the following (for the case of the unit ball, see \cite{Constantin, Gu} and the references therein).
\btheo\label{theo:Carlembed1}
Suppose that $\mu$ is a positive measure on $\mathcal D$, $1\le p\le q<\infty$, $\nu>\frac{n}{r}-1$. Then there
exists a constant $C>0$ such that
\begin{equation}\label{eq:Carlembed11}
\int_{\mathcal D}|f(z)|^qd\mu (z)\le C\|f\|_{p,\nu}^q,\,\,\,f\in A_\nu^p(\mathcal D)
\end{equation}
if and only if
\begin{equation}\label{eq:Carlembed12}
\mu(B_\delta(z))\le C'\Delta^{\frac{q}{p}(\nu+\frac{n}{r})}(\im z),\,\,\,z\in \mathcal D,
\end{equation}
for some $C'>0$ independent of $z$, and for some $\delta\in (0,1)$.
\etheo
\Proof
Let us first prove that (\ref{eq:Carlembed12}) is sufficient. For $f\in \mathcal {H}(\mathcal D)$, we have from Lemma \ref{lem:meanvalue}
that $$|f(z)|^q\le \int_{B_{1/2}(z)}|f(\zeta)|^q\frac{dV(\zeta)}{\Delta^{2n/r}(\im \zeta)}.$$
Using the latter and (\ref{eq:Carlembed12}), we obtain
\Beas
\int_{\mathcal D}|f(z)|^qd\mu (z) &\le& C\int_{\mathcal D}\left(\int_{B_{1/2}(z)}|f(\zeta)|^q\frac{dV(\zeta)}{\Delta^{2n/r}(\im \zeta)}\right)d\mu(z)\\
&\le& C\int_{\mathcal D}\left(C\int_{\mathcal D}\chi_{B_{1/2}(z)}(\zeta)d\mu(z)\right)|f(\zeta)|^q\frac{dV(\zeta)}{\Delta^{2n/r}(\im \zeta)}\\ &\le&
C\|f\|_{p,\nu}^{q-p}\int_{\mathcal D}|f(\zeta)|^p\mu(B_\delta(\zeta))\Delta^{-\frac{q}{p}(\nu+\frac{n}{r})}(\zeta)dV_\nu(\zeta)\\ &\le& C\|f\|_{p,\nu}^q.
\Eeas
Where we used the fact that $\chi_{B_{1/2}(z)}(\zeta)\le \chi_{B_\delta(\zeta)}(z)$ for each $z\in \mathcal D$ and each $\zeta\in \mathcal D$, with $\delta\geq 1/2$ and the pointwise estimate of functions of $A_\nu^p(\mathcal D)$ given by (\ref{eq:pointwiseestimberg}).
\vskip .2cm
Let us now suppose that (\ref{eq:Carlembed11}) holds and prove that in this case, (\ref{eq:Carlembed12}) necessarily holds. For this, we test (\ref{eq:Carlembed11}) with the function $f_\zeta(z)=k_\nu(z,\zeta)^\frac{2}{p}$ where $k_\nu$ is defined in (\ref{eq:normalrepkern}). One easily checks using \cite[Lemma 3.20]{BBGNPR} that
 $f \in A_\nu^p(\mathcal D)$ with norm $\|f_\zeta\|_{p,\nu}=C_{\nu p}.$ Now, for $\da$ sufficiently small, we get in view of Corollary \ref{kor} and (\ref{eq:Carlembed11})
\Beas
V_\nu(B_\delta(\zeta))^{-\frac{q}{p}}\mu(B_\delta(\zeta)) &\simeq& \int_{B_\da(\zeta)}V_\nu(B_\delta(z))^{-\frac{q}{p}}d\mu(z)\\&\leq& C_{pq}(1-c\da)\int_{B_\da(\zeta)}|k_\nu(\zeta,z)|^{2\frac{q}{p}}d\mu(z)\\&\simeq&
\int_{B_\delta(\zeta)}|f_\zeta(z)|^qd\mu(z)\\ &\le& C(p,q,\da)\int_{\mathcal D}|f_\zeta(z)|^qd\mu(z)\\ &\le& C(p,q,\da)\|f_\zeta\|^q_{p,\nu}=C_{\nu p q \da}.
\Eeas
The proof is complete.

\ProofEnd
\begin{definition}
When $p=q$, we call the measures satisfying (\ref{eq:Carlembed12}) Bergman-Carleson measures.
\end{definition}
We next consider the case $1\le q< p<\infty$. We follow the method of Luecking \cite{Luecking1} that uses Kinchine's inequality. We recall that the Rademacher functions $r_n$ are given by $r_n(t)=r_0(2^nt)$, for $n\ge 1$ and $r_0$ is defined as follows
$$
r_0(t)=\left\{ \begin{matrix} 1 &\text{if }& 0\le t-[t]<1/2\\
      -1 & \text{ if } & 1/2\le t-[t]<1
                                  \end{matrix} \right.
$$
$[t]$ is the smallest integer $k$ such $k\le t<k+1$.
\blem[Kinchine's inequality]\label{lem:Kinchine}
For $0 < p <\infty$ there exist constants $0 <L_p\le M_p <\infty$ such that, for all natural numbers $m$ and all complex
numbers $c_1, c_2,\cdots, c_m,$ we have
$$L_p\left(\sum_{j=1}^m|c_j|^2\right)^{p/2}\le \int_0^1\left|\sum_{j=1}^mc_jr_j(t)\right|^pdt\le M_p\left(\sum_{j=1}^m|c_j|^2\right)^{p/2}.$$
\elem

We have the following embedding with loss.
\btheo\label{theo:Carlembed2}
Suppose that $\mu$ is a positive measure on $\mathcal D$,  $\nu>\frac{n}{r}-1$. Then the following assertions
hold.
\begin{itemize}
\item[(i)] Let $1\le q<p <\infty.$ If the function
$$\mathcal {D}\ni z\mapsto \frac{\mu(B_\delta(z))}{\Delta^{\nu+\frac{n}{r}}(\im z)}$$
belongs to $L_\nu^{s}(\mathcal D)$, $s=\frac{p}{p-q}$ for some $\delta\in (0,1)$, then
there
exists a constant $C>0$ such that
\begin{equation}\label{eq:Carlembed21}
\int_{\mathcal D}|f(z)|^qd\mu (z)\le C\|f\|_{p,\nu}^q,\,\,\,f\in A_\nu^p(\mathcal D).
\end{equation}
\item[(ii)] Let $1\le q<p <\infty.$ If for some $\alpha$ large enough, the Bergman projection $P_\alpha$ is bounded on $L_\nu^p(\mathcal D)$ and (\ref{eq:Carlembed21}) holds, then
the function
$$\mathcal {D}\ni z\mapsto \frac{\mu(B_\delta(z))}{\Delta^{\nu+\frac{n}{r}}(\im z)}$$
belongs to $L_\nu^{s}(\mathcal D)$, $s=\frac{p}{p-q}$ for some $\delta\in (0,1)$.
\end{itemize}
\etheo
\Proof
Let us start with the sufficiency. For $\frac{1}{2}\le \delta<1$, we obtain as in the first part of the proof of Theorem \ref{theo:Carlembed1}
that
\Beas \int_{\mathcal D}|f(z)|^qd\mu (z) &\le& C\int_{\mathcal D}|f(\zeta)|^q\mu(B_\delta(\zeta))\frac{dV(\zeta)}{\Delta^{2n/r}(\im \zeta)}\\ &=&
C\int_{\mathcal D}|f(\zeta)|^q\frac{\mu(B_\delta(\zeta))}{\Delta^{\nu+n/r}(\im \zeta)}dV_\nu(\zeta).
\Eeas
It follows from the H\"older inequality and the hypothesis that
\Beas
\int_{\mathcal D}|f(z)|^qd\mu (z) &\le& \left\|\frac{\mu(B_\delta(\zeta))}{\Delta^{\nu+n/r}(\im \zeta)}\right\|_{L_\nu^s}\|f\|_{p,\nu}^q\\ &\le& C \|f\|_{p,\nu}^q
\Eeas
that is (\ref{eq:Carlembed21}) holds.

To prove the converse, we recall with Theorem \ref{theo:atomdecompo} that if the (weighted) Bergman projection $P_\alpha$ ($\alpha$ large enough) is bounded on $L_\nu^p(\mathcal D)$, then given a $\delta$-lattice $\{z_j\}$ (with $\delta$ small enough), every function $f\in A_\nu^p(\mathcal D)$ may be written as
$$f(z)=\sum_j\lambda_jK_\alpha (z,z_j)\Delta^{\alpha+n/r}(\im z_j)$$ where $\{\lambda_j\}\in l_{\nu}^p$ that is
$$\sum_j|\lambda_j|^p\Delta^{\nu+n/r}(\im z_j)<\infty$$ with $$\|f\|_{p,\nu}\le C\left(\sum_j|\lambda_j|^p\Delta^{\nu+n/r}(\im z_j)\right)^{1/p}.$$
It follows that if (\ref{eq:Carlembed21}) holds, then
\Beas \int_{\mathcal D}\left|\sum_j\lambda_jK_\alpha (z,z_j)\Delta^{\alpha+n/r}(\im z_j)\right|^qd\mu(z) &\le& C\|f\|_{p,\nu}^q \le C\|\{\la_j
\}\|_{\ell_\nu^p}^q 
\Eeas
Replacing $\lambda_j$ by $\lambda_jr_j(t)$ yields
\begin{equation}\label{eq:RadfunctionAdd} \int_{\mathcal D}\left|\sum_j\lambda_jr_j(t)K_\alpha (z,z_j)\Delta^{\alpha+n/r}(\im z_j)\right|^qd\mu(z) \le C\|\{\la_jr_j(t)
\}\|_{\ell_\nu^p}^q. 
\end{equation}
By Kinchine's inequality, we have
\Beas \int_{0}^1\left|\sum_j\lambda_jr_j(t)K_\alpha (z,z_j)\Delta^{\alpha+n/r}(\im z_j)\right|^qdt  \ge L_q\left(\sum_j|\lambda_j|^2|K_\alpha (z,z_j)|^2\Delta^{2(\alpha+n/r)}(\im z_j)\right)^{q/2}.
\Eeas
Now, integrating both sides of (\ref{eq:RadfunctionAdd}) from $0$ to $1$ with respect to $dt$, we obtain with the help of Fubini's theorem
that
\Bea\label{kin2}
L_q\int_{\mathcal D}\left(\sum_j|\lambda_j|^2|K_\alpha (z,z_j)|^2\Delta^{2(\alpha+n/r)}(\im z_j)\right)^{q/2}d\mu(z)\leq C\|\{\la_j
\}\|_{\ell_\nu^p}^q .
\Eea
Next we observe that
\Beas
\langle \{|\lambda_j|^q\},\left\{\frac{\mu(B_j)}{\Delta^{\nu+n/r}(\im \zeta_j)}\right\}\rangle_\nu &=& \sum_j|\lambda_j|^q\frac{\mu(B_j)}{\Delta^{\nu+n/r}(\im \zeta_j)}\Delta^{\nu+n/r}(\im \zeta_j)\\ &=&\int_{\mathcal D}\sum_j|\lambda_j|^q\chi_{B_j}(z)d\mu(z).\Eeas
Observe that for $\frac 2q\ge 1,$ by the H\"older inequality, we have using the fact that each point in $\mathcal D$ belongs to at most $N$ balls $B_j$ that for all $z\in\mathcal{D},$
\Beas
\sum_j|\lambda_j|^q\chi_{B_j}(z) &\leq& \left(\sum_j|\lambda_j|^2\chi_{B_j}(z)\right)^\frac{q}{2}\left(\sum_j\chi_{B_j}(z)\right)^\frac{2-q}{2}\leq N^\frac{2-q}{2}\left(\sum_j|\lambda_j|^2\chi_{B_j}(z)\right)^\frac{q}{2},
\Eeas
and that for $\frac 2q< 1$,
$$\sum_j|\lambda_j|^q\chi_{B_j}(z)\le \left(\sum_j|\lambda_j|^2\chi_{B_j}(z)\right)^{\frac q2}.$$
Hence from (\ref{kin2}), we obtain
\Beas
\langle \{|\lambda_j|^q\},\left\{\frac{\mu(B_j)}{\Delta^{\nu+n/r}(\im \zeta_j)}\right\}\rangle_\nu &\le& C\int_{\mathcal D}\left(\sum_j|\lambda_j|^2\chi_{B_j}(z)\right)^{q/2}d\mu(z)\\ &\le& C
\int_{\mathcal D}\left(\sum_j|\lambda_j|^2|K_\alpha (z,z_j)|^2\Delta^{2(\alpha+n/r)}(\im z_j)\right)^{q/2}d\mu(z)\\ &\le& C\|\{\la_j
\}\|_{\ell_\nu^p}^q .
\Eeas
We deduce that the sequence $\{\frac{\mu(B_j)}{\Delta^{\nu+n/r}(\im \zeta_j)}\}$ belongs to $l_\nu^s$ ($s=\frac{p}{p-q}$) which is the dual of
 $l_\nu^{\frac{p}{q}}$ under the pairing $$\langle \{a_j\},\{b_j\}\rangle_\nu:=\sum_ja_j\overline {b_j}\Delta^{\nu+n/r}(\im \zeta_j).$$ That is
$$\sum_j\left(\frac{\mu(B_j)}{\Delta^{\nu+n/r}(\im \zeta_j)}\right)^s\Delta^{\nu+n/r}(\im \zeta_j)<\infty.$$
It follows from Lemma \ref{lem:integraldiscretizationAverBer} that
$$\int_{\mathcal D}\left(\frac{\mu(B_\delta(z))}{\Delta^{\nu+n/r}(\im z)}\right)^sdV_\nu(z)<\infty.$$
The proof is complete.
\ProofEnd
\subsection{Ces\`aro-type operators from $A_\nu^p(\mathcal D)$ to $\mathcal {B}_\nu^q(\mathcal D)$}
We give here a natural application of Carleson embeddings adapted to our setting. We recall that for $1\le p<\infty$, $\nu\in \mathbb {R}$,
$\mathcal {B}_\nu^p(\mathcal D)$ is the set of equivalence classes in $\mathcal {H}_n(\mathcal D)$ of functions
$f\in \mathcal {H}(\mathcal D)$ such that
$$\int_{\mathcal D}|\Delta^n(\im z)\Box^nf(z)|^pdV_\nu (z)<\infty.$$
We have the following consequence of Theorem \ref{theo:Carlembed1} and Theorem \ref{theo:Carlembed2}.
\bcor\label{cor:embedCesaro}
Let $\nu>\frac{n}{r}-1$, $g\in \mathcal {H}(\mathcal D)$. Then the following assertions hold.
\begin{itemize}
\item[(i)] If $1\le p\le q<\infty$, then $T_g$ is bounded from $A_\nu^p(\mathcal D)$ to $\mathcal {B}_\nu^q(\mathcal D)$, if and only if
$$|\Box^ng(\zeta)|\le C\Delta^{(\frac{1}{p}-\frac{1}{q})\nu+(\frac{1}{p}-\frac{1}{q})\frac nr-n}(\im \zeta), \,\,\,\zeta\in \mathcal D.$$
\item[(ii)] If $1\le q<p<\infty$ and the Bergman projection $P_\alpha$ is bounded on $L_\nu^p(\mathcal D)$ for some $\alpha$ large enough, then $T_g$ is bounded
from $A_\nu^p(\mathcal D)$ to $\mathcal {B}_\nu^q(\mathcal D)$ if and only if $g\in \mathcal {B}_\nu^{s}(\mathcal D)$, $\frac{1}{s}=\frac{1}{q}-\frac{1}{p}$.
\end{itemize}
\ecor
\Proof
That $T_g$ is bounded from $A_\nu^p(\mathcal D)$ to $\mathcal {B}_\nu^q(\mathcal D)$ is equivalent in saying that there exists
a positive constant $C$ such that
\begin{equation}\label{eq:Cesarobound}
\int_{\mathcal D}|f(z)|^q|\Box^ng(z)|^q\Delta^{nq+\nu-\frac{n}{r}}(\im z)dV(z)\le C\|f\|_{p,\nu}^q.
\end{equation}
If $1\le p\le q<\infty$, then by Theorem \ref{theo:Carlembed1}, (\ref{eq:Cesarobound}) is equivalent in saying that there exists $C>0$ such that for any
Bergman ball $B_\delta(\zeta)$, $\zeta\in \mathcal D$, $\delta\in (0,1)$,
\begin{equation}\label{eq:Cesaroembedd1}\mu(B_\delta(\zeta))\le C\Delta^{\frac{q}{p}(\nu+n/r)}(\im \zeta)
\end{equation}
where $d\mu(z)=|\Box^ng(z)|^q\Delta^{nq+\nu-\frac{n}{r}}(\im z)dV(z)$. To finish, we only need to prove that the latter is equivalent to
\begin{equation}\label{eq:Cesaroembedd2}|\Box^ng(\zeta)|\le C\Delta^{(\frac{1}{p}-\frac{1}{q})\nu+(\frac{1}{p}-\frac{1}{q})\frac nr-n}(\im \zeta).
\end{equation}

That (\ref{eq:Cesaroembedd1}) implies (\ref{eq:Cesaroembedd2}) is direct from the mean value property. Let us prove the other implication. We obtain using (\ref{eq:Cesaroembedd2}) that
\Beas
\mu(B_\delta(\zeta)) &=& \int_{B_\delta(\zeta)}|\Box^ng(w)|^q\Delta^{nq+\nu-\frac{n}{r}}(\im w)dV(w)\\ &\le& C\int_{B_\delta(\zeta)}\Delta^{\frac{q}{p}(\nu+n/r)-2\frac{n}{r}}(\im w)dV(w)\\ &\le& C\Delta^{\frac{q}{p}(\nu+n/r)}(\im \zeta).
\Eeas
\vskip .2cm
Now if $1\le q<p<\infty$, then by Theorem \ref{theo:Carlembed2}, (\ref{eq:Cesarobound}) is equivalent in saying that for some $\delta\in (0,1)$,
the function $\zeta\mapsto \frac{\mu(B_\delta(\zeta))}{\Delta^{\nu+n/r}(\im \zeta)}$ ($d\mu(z)=|\Box^ng(z)|^q\Delta^{nq+\nu-\frac{n}{r}}(\im z)dV(z)$) belongs to $L_\nu^{\frac{p}{p-q}}(\mathcal D)$. We can once more suppose that $\delta<\frac{1}{2}$. To conclude, we only need to prove that
$$\left\|\frac{\mu(B_\delta(\cdot))}{\Delta^{\nu+n/r}(\im \cdot)}\right\|_{\frac{p}{p-q},\nu}^{\frac{p}{p-q}}\backsimeq \int_{\mathcal D}|\Delta^{n}(\im z)\Box^ng(z)|^{\frac{pq}{p-q}}dV_\nu(z).$$

By the mean value property we get $$\int_{\mathcal D}|\Delta^{n}(\im z)\Box^ng(z)|^{\frac{pq}{p-q}}dV_\nu(z)\le C\left\|\frac{\mu(B_\delta(\cdot))}{\Delta^{\nu+n/r}(\im \cdot)}\right\|_{\frac{p}{p-q},\nu}^{\frac{p}{p-q}}.$$
 Let us prove the converse inequality.
\Beas
 \left\|\frac{\mu(B_\delta(\cdot))}{\Delta^{\nu+n/r}(\im \cdot)}\right\|_{\frac{p}{p-q},\nu}^{\frac{p}{p-q}} &=& \int_{\mathcal D}\left(\frac{\mu(B_\delta(\zeta))}{\Delta^{\nu+n/r}(\im \zeta)}\right)^{\frac{p}{p-q}}dV_\nu(\zeta)\\ &=& \int_{\mathcal D}\left(\frac{1}{\Delta^{\nu+n/r}(\im \zeta)}\int_{B_\delta(\zeta)}|\Delta^{n}(\im z)\Box^ng(z)|^qdV_\nu(z)\right)^{\frac{p}{p-q}}dV_\nu(\zeta)\\ &\le& \sum_j\int_{B_j}\left(\int_{B_\delta(\zeta)}|\Delta^{n}(\im z)\Box^ng(z)|^qdV_\nu(z)\right)^{\frac{p}{p-q}}\Delta^{-\frac{q}{p-q}(\nu+n/r)}(\im \zeta)d\lambda (\zeta).
 \Eeas

By H\"older's inequality, we have
\Beas\int_{B_\delta(\zeta)}|\Delta^{n}(\im z)\Box^ng(z)|^qdV_\nu(z)&\leq& \left(\int_{B_\delta(\zeta)}|\Delta^{n}(\im z)\Box^ng(z)|^{\frac{pq}{p-q}}dV_\nu(z)\right)^{\frac{p-q}{p}}\left(V_\nu(B_\delta(\zeta))\right)^{\frac{q}{p}}\\&\lesssim &\left(\int_{B_\delta(\zeta)}|\Delta^{n}(\im z)\Box^ng(z)|^{\frac{pq}{p-q}}dV_\nu(z)\right)^{\frac{p-q}{p}}\Delta^{\frac{q}{p}(\nu+n/r)}(\im \zeta).
\Eeas
Hence
 \Beas
   \left\|\frac{\mu(B_\delta(\cdot))}{\Delta^{\nu+n/r}(\im \cdot)}\right\|_{\frac{p}{p-q},\nu}^{\frac{p}{p-q}}&\lesssim & \sum_j\int_{B_j}\left(\int_{B_\delta(\zeta)}|\Delta^{n}(\im z)\Box^ng(z)|^{\frac{pq}{p-q}}dV_\nu(z)\right)d\lambda (\zeta)\\
   &\lesssim & \sum_j\int_{B_j}\left(\int_{B_{2\delta}(z_j)}|\Delta^{n}(\im z)\Box^ng(z)|^{\frac{pq}{p-q}}dV_\nu(z)\right)d\lambda (\zeta)\\
   &\lesssim & \sum_j\int_{B_{2\delta}(z_j)}|\Delta^{n}(\im z)\Box^ng(z)|^{\frac{pq}{p-q}}dV_\nu(z)\\
   &\lesssim & M\int_{\mathcal D}|\Delta^{n}(\im z)\Box^ng(z)|^{\frac{pq}{p-q}}dV_\nu(z);
 \Eeas
 the second inequality above is due to the fact that for $\zeta\in B_j=B_\da(z_j)$ we have $B_\da(\zeta)\subset B_{2\da}(z_j).$
%
The proof is complete.
\ProofEnd
\begin{Remark}
Let us observe that in the case $q<p$, even if we suppose that $p$ and $q$ belong to an interval on which $P_\nu$ extends as bounded operator,
the above result tell us that the range of symbols for which the Ces\`aro-type operators $T_g$ are bounded that is $\mathcal {B}_\nu^s(\mathcal D)$ ($s=\frac{pq}{p-q}$) is in general bigger than the range $A_\nu^s(\mathcal D)$ obtained in dimension one (see \cite{AlemanSis}).
\end{Remark}

Let us finish this section with the following result which does not use the previous Carleson embeddings.
\btheo
Let $1\le p<\infty$, $\nu>\frac{n}{r}-1$. For $g\in \mathcal {H}(\mathcal D)$, $T_g$ is bounded from
$A_\nu^p(\mathcal D)$ to $\mathcal {B}^\infty(\mathcal D)$ if and only if
\begin{equation}\label{eq:CesaroBlock}
\sup_{\zeta\in \mathcal D}\Delta^{n-\frac{1}{p}(\nu+n/r)}(\im \zeta)|\Box^ng(\zeta)|<\infty.
\end{equation}
\etheo
\Proof
Let us prove that (\ref{eq:CesaroBlock}) is sufficient. We recall that there is a constant $C>0$ such that the following pointwise estimate holds
$$|f(z)|\le C\Da^{-\frac{1}{p}(\nu+n/r)}(\im z)\|f\|_{p,\nu},\,\,\,f\in A_\nu^p(\mathcal D).$$ It follows easily that for any $f$ in the
unit ball of $A_\nu^p(\mathcal D)$ and any $\zeta\in \mathcal D$,
\Beas
\Delta^{n}(\im \zeta)|\Box^n(T_gf)(\zeta)| &=& \Delta^{n}(\im \zeta)|f(\zeta)||\Box^ng(\zeta)|\\ &\le& C\Delta^{n-\frac{1}{p}(\nu+n/r)}(\im \zeta)|\Box^ng(\zeta)|<\infty.
\Eeas
That is $T_g$ is bounded from $A_\nu^p(\mathcal D)$ to $\mathcal {B}^\infty$.

\vskip .2cm
Conversely, if $T_g$ is bounded from $A_\nu^p(\mathcal D)$ to $\mathcal {B}^\infty$, then for any $f\in A_\nu^p(\mathcal D)$,
$$\Delta^{n}(\im \zeta)|f(\zeta)||\Box^ng(\zeta)|<\infty, \,\,\,\textrm{for any}\,\,\,\zeta\in \mathcal D.$$
Taking $f(\zeta)=f_z(\zeta)=\Da^{-\frac{2}{p}(\nu+n/r)}(\zeta-\overline {z})\Da^{\frac{1}{p}(\nu+n/r)}(\im z)$ in the later inequality and, then putting in particular $\zeta=z$,  we conclude that
$$\sup_{\zeta\in \mathcal D}\Delta^{n-\frac{1}{p}(\nu+n/r)}(\im \zeta)|\Box^ng(\zeta)|<\infty.$$
The proof is complete.
\ProofEnd

\section{Schatten class membership of Toeplitz and Ces\`aro-type operators}
The aim of this section is to give criteria for Schatten class
membership of Toeplitz and Ces\`aro-type operators on the weighted Bergman space
$A_\nu^2(\mathcal D)$.

Let  $\mathcal{H}$ be a Hilbert
space and $\|\,\|$ its associated norm. Let $\mathcal{B}(\mathcal{H})$ and
$\mathcal{K}(\mathcal{H})$ denote the
spaces of bounded and compact linear operators on
$\mathcal{H}$ respectively. It is
well known that any operator $T\in
\mathcal{K}(\mathcal{H})$ has a Schmidt
decomposition, that is there exist an orthonormal base
$\{e_j\}$ of $\mathcal{H}$ and a sequence $\{\lambda_j\}$
of complex numbers with $\lambda_j\rightarrow 0$, such that
\begin{equation}
Tf=\sum_{j=0}^{\infty}\lambda_j\langle f,e_j\rangle e_j,\,\,\,\,\,\,\,\,\,\,
f\in \mathcal{H} .\end{equation}

For $1\le p <\infty$, a
compact operator $T$ with such decomposition belongs to
the Schatten-Von Neumann p-class
$\mathcal{S}_p(\mathcal{H})$, if and only
if
$$||T||_{\mathcal{S}_p}:=(\sum_{j=0}^{\infty}|\lambda_j|^p)^{\frac{1}{p}}<\infty.$$


For $p=1$, $\mathcal{S}_1=\mathcal{S}_1(\mathcal{H})$ is
the trace class and for $T\in \mathcal{S}_1$, the trace of
$T$ is defined by $$Tr(T)=\sum_{j=0}^{\infty}\langle Te_j,e_j\rangle$$
where $\{e_j\}$ is any orthonormal basis of the Hilbert
space $\mathcal{H}$. We will denote by $\mathcal
{S}_\infty(\mathcal H)$ the set of all bounded linear
operators on $\mathcal H$.

We recall that a compact operator T on $\mathcal H$ belongs
to the Schatten class $\mathcal {S}_p$ if the positive operator $ (T^*T)^{1/2}$ belongs
to $\mathcal {S}_p$, where $T^*$ denotes the adjoint of $T$. In this case, we have $||T||_{\mathcal{S}_p}=||(T^*T)^{1/2}||_{\mathcal{S}_p}$.
\vskip .2cm
We shall use also the fact that if $T$ is a compact operator on $\mathcal H$, and $p\ge 2$, then that $T\in \mathcal S_p$ is equivalent to
$$\sum_{j}\|Te_j\|^p<\infty$$ for any orthonormal set $\{e_j\}$ in $\mathcal H$ (see \cite{Zhu1}).
\vskip .2cm
Let us start with an elementary result which gives that for a positive Borel measure $\mu$ on $\mathcal D$, the Toeplitz operator $T_\mu$
belongs to $\mathcal {S}_\infty(A_\nu^2(\mathcal D))$, $\nu>\frac{n}{r}-1$.
\blem\label{lem:boundedToeplitz}
Let $\mu$ be a positive Borel measure on $\mathcal D$ and $\nu>\frac{n}{r}-1$. Then the Toeplitz operator $T_\mu$ is bounded on $A_\nu^2(\mathcal D)$ if and only $\mu$ is a Bergman-Carleson measure.
\elem
\Proof
The proof is pretty classical, we give it here for completeness. First, we suppose that $T_\mu$ is bounded on $A_\nu^2(\mathcal D)$. Then for any
$f\in A_\nu^2(\mathcal D)$, using Fubini's theorem and reproducing properties of the Bergman kernel, we have that
$$\int_{\mathcal D}|f(z)|^2d\mu(z)=\left| \langle T_\mu f,f\rangle_\nu\right|\le C\|f\|_{2,\nu}^2.$$
That is $\mu$ is a Bergman-Carleson measure.
\vskip .2cm
Conversely, if $\mu$ is a Bergman-Carleson measure, then using Theorem \ref{theo:Carlembed1} we obtain that for any $f,g\in A_\nu^2(\mathcal D)$,
$$\left|\langle T_\mu f,g\rangle_\nu\right|\le \int_{\mathcal D}|f(z)g(z)|d\mu(z)\le C\|fg\|_{1,\nu}\le C\|f\|_{2,\nu}\|g\|_{2,\nu}.$$
It follows that $T_\mu$ is bounded on $A_\nu^2(\mathcal D)$.
\ProofEnd

\vskip .2cm
We recall that the Berezin transform of a measure $\mu$ is the function defined on $\mathcal D$ by
$$\tilde {\mu}(z):=\int_{\mathcal D}|k_\nu(z,\zeta)|^2d\mu(\zeta)$$
where the function $k_\nu$ is introduced in (\ref{eq:normalrepkern}). Observe that
\begin{equation}\label{eq:Berezininnerprod}\tilde {\mu}(z)=\langle T_\mu k_\nu(\cdot,z),k_\nu(\cdot,z)\rangle_\nu.
\end{equation}
\vskip .2cm
Let us prove the following criteria of Schatten classes membership of the Toeplitz operators.

\btheo\label{theo:SchattenToep}
Let $\mu$ be a positive Borel measure on $\mathcal D$, and $\nu>\frac{n}{r}-1$. Then for $p\ge 1$, the following assertions are equivalent
\begin{itemize}
\item[(i)] The Toeplitz operator $T_\mu$ belongs
to the Schatten class $\mathcal {S}_p(A_\nu^2(\mathcal D))$.
\item[(ii)] For any $\delta$-lattice ($\delta\in (0,1)$) $\{\zeta_j\}_{j\in \mathbb N}$ in the Bergman metric of $\mathcal D$, the sequence $\{\hat {\mu}_\delta(\zeta_j)\}$ belongs to $l^p$, that is \begin{equation}\label{eq:SchattenToepcond}
\sum_j\left(\frac{\mu(B_j)}{\Da^{\nu+n/r}(\im \zeta_j)}\right)^p<\infty.
\end{equation}
\item[(iii)] For any $\beta\in (0,1)$, the function $z\mapsto \hat {\mu}_\beta(z)$ belongs to $L^p(\mathcal {D},d\lambda)$.
\item[(iv)] $\tilde {\mu}\in L^p(\mathcal {D},d\lambda)$.

\end{itemize}
\etheo
\Proof
The equivalences $(ii)\Leftrightarrow (iii)\Leftrightarrow (iv)$ are given by Lemma \ref{lem:integraldiscretizationAverBer}. We prove that $(i)\Leftrightarrow (ii)$.

$(i)\Rightarrow (ii)$: Clearly, using Corollary \ref{kor} and (\ref{eq:Berezininnerprod}), we obtain
\Beas
\sum_{j}\left(\frac{\mu(B_j)}{\Da^{\nu+n/r}(\im \zeta_j)}\right)^p &\le& C\sum_j\left(\int_{B_j}|k_\nu(z,\zeta_j)|^2d\mu(z)\right)^p\\ &\le& C\sum_j\left(\int_{\mathcal D}|k_\nu(z,\zeta_j)|^2d\mu(z)\right)^p\\ &=& C\sum_j\left|\langle T_\mu k_\nu(\cdot,\zeta_j),k_\nu(\cdot,\zeta_j)\rangle\right|^p.
\Eeas
To conclude, it is enough
to prove that $k_\nu(\cdot,\zeta_j)$ is the image of
an orthonormal sequence $e_j$ through a bounded linear operator on $L_\nu^2(\mathcal D)$.

Consider
$$ e_j(z)=C_{n,\nu}\Da^{-\frac{1}{2}(\nu+n/r)}(\im z)\chi_{B_j'}(z).$$
Clearly, $C_{n,\nu}P_{\nu} e_j=k_\nu(\cdot,\zeta_j)$ and $\|e_j\|_{2,\nu}=1$ for an appropriate choice of $C_{n,\nu}$.

\vskip .2cm
$(ii)\Rightarrow (i)$: We proceed by interpolation. We observe that from
Proposition \ref{prop:Carlembedlattice}, we clearly have that if the sequence $\{\frac{\mu(B_j)}{\Da^{\nu+n/r}(\im \zeta_j)}\}_{j\in \N}$ belongs
to $l^\infty$, that is, $$\sup_j\frac{\mu(B_j)}{\Da^{\nu+n/r}(\im \zeta_j)}\le C<\infty,$$ then $T_\mu\in \mathcal {S}_\infty(A_\nu^2(\mathcal D))$. To conclude, we only need to prove that if the sequence $\{\frac{\mu(B_j)}{\Da^{\nu+n/r}(\im \zeta_j)}\}_{j\in \N}$ is in $l^1$ i.e.
$$\sum_j\frac{\mu(B_j)}{\Da^{\nu+n/r}(\im \zeta_j)}<\infty,$$
then $T_\mu$ belongs to $\mathcal {S}_1(A_\nu^2(\mathcal D))$. For this, we only need to check that for any orthonormal basis
$\{e_k\}_{k\in \N}$ of $A_\nu^2(\mathcal D)$, $$\sum_{k\in \N}|\langle T_\mu e_k,e_k\rangle_\nu|<\infty.$$

Clearly,
\Beas
\sum_{k\in \mathbb N}\langle T_\mu e_k,e_k\rangle_\nu &=& \sum_{k\in \N}\int_{\mathcal D}|e_k(z)|^2d\mu(z)\\ &=& \int_{\mathcal D}K_\nu(z,z)d\mu(z)\\ &\le&
C\sum_{j}\frac{\mu(B_j)}{\Da^{\nu+n/r}(\im \zeta_j)}<\infty.
\Eeas
Hence, $$T_\mu\in \mathcal {S}_p(A_\nu^2(\mathcal D)), p\ge 1, \,\,\, \textrm{if}\,\,\, \sum_j\left(\frac{\mu(B_j)}{\Da^{\nu+n/r}(\im \zeta_j)}\right)^p<\infty.$$
The proof is complete.
\ProofEnd

Let us recall that $\mathcal {B}^p(\mathcal D)$ is the subset of $\mathcal {H}_n(\mathcal D)$ consisting of functions $f$ such that
$\Da^n\Box^nf\in L^p(\mathcal {D}, d\lambda)=L^p(\mathcal {D}, \frac{dV(z)}{\Da^{2\frac{n}{r}}(\im z)})$.

We have the following result.
\btheo\label{theo:SchattenCesaro}
Let $p\ge 2$, $\nu>\frac{n}{r}-1$. If $g$ is a given holomorphic function in $\mathcal D$, then the Ces\`aro-type operator $T_g$ belongs to
$\mathcal {S}_p(A_\nu^2(\mathcal D))$ if and only if $g\in \mathcal {B}^p(\mathcal D)$.
\etheo
\Proof
Let us first suppose that $g$ is such that $T_g\in \mathcal {S}_p(A_\nu^2(\mathcal D))$. It follows that for any orthonormal basis $\{e_j\}_{j\in \N}$
of $A_\nu^2(\mathcal D)$,
$$\sum_j \|T_ge_j\|_{2,\nu}^p<\infty.$$
But we have
\Beas \|T_ge_j\|_{2,\nu}^2 &\simeq& \|\Da^n\Box^n\left(T_ge_j\right)\|_{2,\nu}^2\\ &=&\int_{\mathcal D}|e_j(z)|^2|\Box^ng(z)|^2\Da^{2n+\nu-n/r}(\im z)dV(z)\\
&=& \int_{\mathcal D}|e_j(z)|^2d\mu(z)\\ &=& |\langle T_\mu e_j,e_j\rangle_\nu|,
\Eeas
where $d\mu(z)=|\Box^ng(z)|^2\Da^{2n+\nu-n/r}(\im z)dV(z)$. Thus,
$$\sum_j|\langle T_\mu e_j,e_j\rangle_\nu|^{p/2}=\sum_j \|T_ge_j\|_{2,\nu}^p<\infty$$
for any orthonormal basis of $A_\nu^2(\mathcal D)$. This is equivalent in saying that the Toeplitz operator $T_\mu$ belongs to $\mathcal {S}_{p/2}(A_\nu^2(\mathcal D))$ which by Theorem \ref{theo:SchattenToep} is equivalent to the condition
$$\sum_j\left(\frac{\mu(B_j)}{\Da^{\nu+n/r}(\im \zeta_j)}\right)^{p/2}<\infty$$
for some $\delta$-lattice $\{\zeta_j\}_{j\in \N}$ with $\delta\in (0,1)$, $B_j=B_\delta(\zeta_j)$.
Hence, to conclude, we only need to prove that the latter implies that $g\in \mathcal {B}^p(\mathcal D)$.

Using Lemma \ref{5.3}, we obtain
\Beas
\int_{\mathcal D}|\Da^n(\im z)\Box^ng(z)|^pd\lambda(z) &\backsimeq& 
\sum_j \left(|\Box^ng(\zeta_j)|^2\Da^{2n}(\im \zeta_j)\right)^{p/2}\\ &\lesssim&
\sum_j\left(\int_{B_j}|\Box^ng(z)|^2\Da^{2n}(\im z)\frac{dV(z)}{\Da^{2n/r}(\im z)}\right)^{p/2}\\ &\simeq& \sum_j\left(\frac{1}{\Da^{\nu+n/r}(\im \zeta_j)}\int_{B_j}d\mu(z)\right)^{p/2}\\ &=& \sum_j\left(\frac{\mu(B_j)}{\Da^{\nu+n/r}(\im \zeta_j)}\right)^{p/2}\\ &<& \infty.
\Eeas

We next prove that the condition $g\in \mathcal {B}^p(\mathcal D)$ is sufficient for $T_g$ to belong to $\mathcal {S}_p(A_\nu^2(\mathcal D))$, $p\ge 2$. Once more, we proceed by interpolation. We first observe that by Corollary \ref{cor:embedCesaro}, we have that if $g\in \mathcal {B}^\infty(\mathcal D)$, then $T_g\in \mathcal {S}_\infty(A_\nu^2(\mathcal D))$. Thus, to finish the proof, we only need to prove that if $g\in \mathcal {B}^2(\mathcal D)$, then $T_g\in \mathcal {S}_2(A_\nu^2(\mathcal D))$.

\vskip .2cm
Let $\{e_j\}_{j\in \N}$ be an orthonormal basis of $A_\nu^2(\mathcal D)$. Then
\Beas
\sum_j \|T_ge_j\|_{2,\nu}^2 &\simeq& \sum_j\|\Da^n\Box^n\left(T_ge_j\right)\|_{2,\nu}^2\\ &=&\sum_j\int_{\mathcal D}|e_j(z)|^2|\Box^ng(z)|^2\Da^{2n+\nu-n/r}(\im z)dV(z)\\ &=& \int_{\mathcal D}K_\nu(z,z)|\Box^ng(z)|^2\Da^{2n+\nu-n/r}(\im z)dV(z)\\ &=& \int_{\mathcal D}|\Da^n(\im z)\Box^ng(z)|^2\frac{dV(z)}{\Da^{2n/r}(\im z)}<\infty.
\Eeas
The proof is complete.

\ProofEnd


\bibliographystyle{plain}

\end{document}